\documentclass[11pt]{amsart}
\usepackage{amssymb}
\newcommand{\cc}{\mathbb{C}}
\newcommand{\zz}{\mathbb{Z}}
\newcommand{\dbar}{\overline{\partial}} 
\DeclareMathOperator{\dist}{dist}

\newtheorem{theorem}{Theorem}
\newtheorem{lemma}[theorem]{Lemma}
\newtheorem{definition}[theorem]{Definition}

\title[$\dbar$-Integration in homogeneous varieties] 
{An explicit $\dbar$-integration formula for weighted 
homogeneous varieties}

\author{J. Ruppenthal}
\author{E. S. Zeron}

\address{Department of Mathematics, University of Michigan, 
530 Curch Street, Ann Arbor, MI 48109, USA.}
\email{jean@math.uni-bonn.de}

\address{Depto. Matem\'aticas, CINVESTAV, Apartado 
Postal 14-740, M\'exico D.F., 07000, M\'exico.}
\email{eszeron@math.cinvestav.mx}

\date{January 10, 2008}
\thanks{
The first author was supported by a fellowship within the 
Postdoc-Programme of the German Academic Exchange 
Service (DAAD). The second author was supported by 
Cinvestav(Mexico) and Conacyt-SNI(Mexico)}
\subjclass[2000]{32F20, 32W05, 35N15}
\keywords{Cauchy-Riemann equations, H\"older estimates, 
$L^2$-estimates, resolution of singularities}

\begin{document}
\begin{abstract}
Let $\Sigma$ be a weighted homogeneous (singular) subvariety of $\cc^n$. 
The main objective of this paper is to present an explicit formula for 
solving the $\dbar$-equation $\lambda=\dbar{g}$ on the regular part of 
$\Sigma$, where $\lambda$ is a $\dbar$-closed $(0,1)$-form with compact 
support. This formula will then be used to give H\"older estimates for 
the solution in case $\Sigma$ is homogeneous (a cone) with an isolated 
singularity. Finally, a slight modification of our formula also gives an 
$L^2$-bounded solution operator in case $\Sigma$ is pure $d$-dimensional 
and homogeneous.
\end{abstract}

\maketitle
\section{Introduction}

As it is well known, solving the $\dbar$-equation forms a main part of 
complex analysis, but also has deep consequences on algebraic geometry, 
partial differential equations and other areas. In general, it is not easy 
to solve the $\dbar$-equation. The existence of solutions depends mainly 
on the geometry of the variety on which the equation is considered. There 
is a vast literature about this subject on smooth manifolds, both in books 
and papers \cite{HL,LM}, but the theory on singular varieties has been 
developed only recently.

Let $\Sigma$ be a singular subvariety of the space $\cc^n$, and $\lambda$ 
be a $\dbar$-closed differential form well defined on the regular part of 
$\Sigma$. Forn{\ae}ss, Gavosto and Ruppenthal have produced a general 
technique for solving the $\dbar$-equation $\lambda=\dbar{g}$ on the 
regular part of $\Sigma$; see \cite{Ga,FG} and \cite{Ru}. On the other 
hand, Acosta, Sol\'{\i}s and Zeron have proposed an alternative technique 
for solving the $\dbar$-equation when $\Sigma$ is a quotient variety; see 
\cite{AZ1,AZ2} and \cite{SZ}.

In both cases, the main strategy is to transfer the problem to some 
non-singular complex space, to solve the $\dbar$-equation in this 
well-known situation, and to carry over the solution to the singular 
variety. The main objective of this paper is to present and analyze an 
explicit formula for calculating solutions $g$ to the $\dbar$-equation 
$\lambda=\dbar{g}$ on the regular part of the original variety $\Sigma$, 
where $\Sigma$ is a weighted homogeneous variety and $\lambda$ is a 
$\dbar$-closed $(0{,}1)$-differential form with compact support. We 
analyze the weighted homogeneous varieties, for they are a main model for 
classifying the singular subvarieties of $\cc^n$. A detailed analysis of 
the weighted homogeneous varieties is done in Chapter~2--\S4 and Appendix 
B of \cite{Di}.

\begin{definition}\label{whs}

Let $\beta\in\zz^n$ be a fixed integer vector with strictly positive 
entries $\beta_k\geq1$. A polynomial $Q(z)$ holomorphic on $\cc^n$ is said 
to be \textbf{weighted homogeneous} of degree $d\geq1$ with respect to 
$\beta$ if the following equality holds for all $s\in\cc$ and $z\in\cc^n$:
\begin{eqnarray}
\label{homo}&&Q(s^{\beta}*z)\,=\,
s^d\,Q(z),\quad\hbox{with~the~action:}\\
\label{action}&&s^{\beta}*(z_1,z_2,...,z_n)\,:=\,
(s^{\beta_1}z_1,s^{\beta_2}z_2,...,s^{\beta_n}z_n).
\end{eqnarray}

Besides, an algebraic subvariety $\Sigma$ in $\cc^n$ is said to be 
\textbf{weighted homogeneous} with respect to $\beta$ whenever $\Sigma$ 
is the zero locus of a finite number of weighted homogeneous polynomials 
$Q_k(z)$ of (maybe different) degrees $d_k\geq1$, but all of them with 
respect to the same fixed vector $\beta$.
\end{definition}

Let $\Sigma\subset\cc^n$ be any subvariety. We use the following notation 
along this paper. The regular part $\Sigma^*=\Sigma_{reg}$ is the complex 
manifold composed by all the regular points of $\Sigma$, and it is always 
endowed with the induced metric; so that $\Sigma^*$ is a Hermitian 
submanifold in $\cc^n$ with corresponding volume element $dV_\Sigma$ and 
induced norm $|\cdot|_\Sigma$ on the Grassmannian $\Lambda{T}^*\Sigma^*$. 
Thus, any Borel-measurable $(0,1)$-form $\lambda$ on $\Sigma^*$ admits a 
representation $\lambda=\sum_kf_kd\overline{z_k}$, where the coefficients 
$f_k$ are Borel-measurable functions on $\Sigma^*$ which satisfy the 
inequality $|f_k(w)|\leq |\lambda(w)|_\Sigma$ for all points 
$w\in\Sigma^*$ and indexes $1\leq{k}\leq{n}$. Notice that such a 
representation is by no means unique. We refer to Lemma 2.2.1 in 
\cite{Ru} for a more detailed treatment of that point. We also 
introduce the $L^2$-norm of a measurable $(p,q)$-form $\aleph$ 
on an open set $U\subset\Sigma^*$ via the formula:
\begin{eqnarray*}
\|\aleph\|_{L^2_{p,q}(U)} &:=&\bigg(\int_{U} 
|\aleph|^2_\Sigma\,dV_\Sigma\bigg)^{1/2}.
\end{eqnarray*}

We can now present the main result of this paper. We assume that the 
$\dbar$-differentials are calculated in the sense of distributions, 
for we work with Borel-measurable functions.

\begin{theorem}[\textbf{Main}]\label{main}
Let $\Sigma$ be a weighted homogeneous subvariety of $\cc^n$ with 
respect to a given vector $\beta\in\zz^n$, where $n\geq2$ and all 
entries $\beta_k\geq1$. Consider a $(0{,}1)$-form $\lambda$ given 
by $\sum_kf_kd\overline{z_k}$, where the coefficients $f_k$ are all 
Borel-measu\-rable functions in $\Sigma$, and $z_1,..., z_n$ are 
the Cartesian coordinates of $\cc^n$. The following function is well 
defined for all $z\in\Sigma$ whenever the form $\lambda$ is bounded 
and has compact support in $\Sigma$:
\begin{equation}\label{fng-1}
g(z):=\sum_{k=1}^n\frac{\beta_k}{2\pi{i}}\int_{w\in\cc}
f_k(w^{\beta}*z)\frac{(\overline{w^{\beta_k}z_k})\,dw
\wedge{}d\overline{w}}{\overline{w}\;(w-1)}.
\end{equation}
Besides, the function $g$ is a solution of the $\dbar$-equation 
$\lambda=\dbar{g}$ on the regular part of $\Sigma$, whenever 
$\lambda$ is also $\dbar$-closed on the regular part of $\Sigma$. 
\end{theorem}

Notice that $g(0)=0$ and that~(\ref{fng-1}) can also be rewritten 
as follows after the change of variables $u=ws$, given $s\in\cc$ 
and $z\in\Sigma$,
\begin{equation}\label{fng-2}
g(s^{\beta}*z)=\sum_{k=1}^n\frac{\beta_k}{2\pi{i}}\int_{u\in\cc}
f_k(u^{\beta}*z)\frac{(\overline{u^{\beta_k}z_k})\,du\wedge
d\overline{u}}{\overline{u}\;(u-s)}.
\end{equation}

We shall prove Theorem \ref{main} in Section \ref{section:main} of this 
paper. Moreover, recalling some main principles of the proof, we also 
deduce anisotropic H\"older estimates for the $\dbar$-equation in the case 
where $\Sigma$ is a homogeneous variety with an isolated singularity at 
the origin. We obviously need to specify the metric on $\Sigma$: Given a 
pair of points $z$ and $w$ in $\Sigma$, we define $\dist_\Sigma(z,w)$ to 
be the infimum of the length of all piecewise smooth curves connecting 
$z$ and $w$ inside $\Sigma$. It is clear that such curves exist in this 
situation, and that the length of each curve can be measured in the 
regular part $\Sigma^*$ or the ambient space $\cc^n$, but both measures 
coincide, for $\Sigma^*$ carries the induced norm. The main result of 
the section \ref{section:hoelder} is the following estimate:

\begin{theorem}[\textbf{H\"older}]\label{hoelder}
In the situation of Theorem \ref{main}, suppose that $\Sigma$ is 
homogeneous (a cone) and has got only one isolated singularity at the 
origin of $\cc^n$, so that each entry $\beta_k=1$ in Definition~\ref{whs}. 
Moreover, assume that the support of the form $\lambda$ is contained in a 
ball $B_R$ of radius $R>0$ and center at the origin. Then, for each 
parameter $0<\theta<1$, there exists a strictly positive constant 
$C_\Sigma(R,\theta)$ which does not depend on $\lambda$ such that the 
following inequality holds for the function $g$ given in (\ref{fng-1}) 
and all points $z$ and $w$ in the intersection $B_R\cap\Sigma$,
\begin{equation}\label{hoelder2}
|g(z)-g(w)|\leq{C}_\Sigma(R,\theta)\cdot 
\dist_\Sigma(z,w)^\theta\cdot\|\lambda\|_\infty.
\end{equation}
\end{theorem}

The notation $\|\lambda\|_\infty$ stands for the essential supremum 
of $|\lambda(w)|_\Sigma$ on $\Sigma$, recall that $\lambda$ is 
bounded and has compact support. Theorem~\ref{hoelder} is proved in 
Section~\ref{section:hoelder}. Finally, similar techniques and a slight 
modification of equation~(\ref{fng-1}), can also be used for producing 
a $\dbar$-solution operator with $L^2$-estimates on homogeneous 
subvarieties with an isolated singularity at the origin.

\begin{theorem}[\textbf{$L^2$-Estimates}]\label{thm:l2} 
Let $\Sigma$ be a pure $d$-dimensional homogeneous (cone) 
subvariety of $\cc^n$, where $n\geq2$ and each entry $\beta_k=1$ in 
Definition~\ref{whs}. Consider a $(0{,}1)$-form $\lambda$ given by 
$\sum_kf_kd\overline{z_k}$, where the coefficients $f_k$ are all square 
integrable functions in $\Sigma$, and $z_1, ..., z_n$ are the Cartesian 
coordinates of $\cc^n$. The following function is well defined for 
almost all $z\in\Sigma$ whenever the form $\lambda$ has compact 
support on $\Sigma$:
\begin{equation}\label{eq:l2}
g(z):=\sum_{k=1}^n\frac{1}{2\pi{i}}\int_{w\in\cc}
f_k(wz)\frac{w^{d-1}\,\overline{z_k}\,dw\wedge 
d\overline{w}}{w-1}.
\end{equation}
The function $g$ is a solution of the $\dbar$-equation $\lambda=\dbar{g}$ 
on the regular part of $\Sigma$, whenever $\lambda$ is also $\dbar$-closed 
on the regular part or $\Sigma$. Finally, assuming that the support of 
$\lambda$ is contained in an open ball $B_R$ of radius $R>0$ and center 
in the origin, there exists a strictly positive constant $C_\Sigma(R,2)$ 
which does not depend on $\lambda$ and such that:
\begin{equation}\label{eq:l3}
\|g\|_{L^2(\Sigma\cap{B_R})}\leq
C_\Sigma(R,2)\cdot\|\lambda\|_{L^2_{0,1}(\Sigma)}.
\end{equation} 
\end{theorem}

We prove this theorem in Section~\ref{section:L2} of this paper. The 
obstructions to solving the $\dbar$-equation with $L^2$-estimates on 
singular complex spaces are not completely understood in general. An 
$L^2$-solution operator is only known for the case when $\Sigma$ 
is a complete intersection of pure dimension $\geq 3$ with isolated singularities only.
This operator was built by Forn{\ae}ss, {\O}vrelid 
and Vassiliadou in \cite{FOV2}, via an extension theorem for the 
$\dbar$-cohomology groups originally presented by Scheja \cite{Sch}. The 
$L^2$-results usually come with some obstructions to the solvability of 
the $\dbar$-equation. For example, different situations are analyzed in 
the works of Diederich, Forn{\ae}ss, {\O}vrelid and Vassiliadou; it is 
shown there that the $\dbar$-equation is solvable with $L^2$-estimates 
for all forms lying in a closed subspace of finite codimension of the 
vector space of all the $\dbar$-closed $L^2$-forms \cite{DFV,Fo,FOV2,OV}.  
Besides, in the paper \cite{FOV1}, the $\dbar$-equation is solved locally 
with some weighted $L^2$-estimates for forms which vanish to a 
sufficiently high order on the singular set of the given varieties.

% Theorem \ref{thm:l2} also demands strong conditions, for only
% $\dbar$-closed $(0,1)$-forms with compact support are considered 
% there. Anyhow, this result can be used to construct an $L^2$-parametrix 
% (a solution operator modulo a compact defect), which would lead again 
% to the solutions on a subspace of finite codimension of the space of 
% all the $\dbar$-closed $L^2$-forms. I FORGOT THE IDEA BEHIND THAT. 
% DO NOT WANT TO CLAIM IT WITHOUT BEING SURE

On the other hand, we propose in Section~\ref{final} of this 
paper a useful technique for generalizing the estimates given in 
Theorems~\ref{hoelder} and~\ref{thm:l2}, so as to consider weighted 
homogeneous subvarieties instead of homogeneous ones.

\section{Proof of Main Theorem}\label{section:main}

Let $\{Q_k\}$ be the set of polynomials on $\cc^n$ which defines the 
algebraic variety $\Sigma$ as its zero locus. The definition of weighted 
homogeneous varieties implies that the polynomials $Q_k(z)$ are all 
weighted homogeneous with respect to the same fixed vector $\beta$. 
Equation~(\ref{homo}) automatically yields that every point $s^{\beta}*z$ 
lies in $\Sigma$ for all $s\in\cc$ and $z\in\Sigma$, and so each 
coefficient $f_k(\cdot)$ in equations~(\ref{fng-1}) and~(\ref{fng-2}) 
is well evaluated in $\Sigma$. Moreover, fixing any point $z\in\Sigma$, 
the given hypotheses imply that the following Borel-measu\-rable 
functions are all bounded and have compact support in $\cc$,
$$w\,\mapsto\,f_k(w^{\beta}*z).$$

Hence, the function $g(z)$ in~(\ref{fng-1}) is well defined for every 
$z\in\Sigma$. Notice that $g(0)=0$, in particular. We shall prove that 
$g(z)$ is also a solution of the equation $\dbar{g}=\lambda$, when the 
$(0{,}1)$-form $\lambda$ is $\dbar$-closed. We may suppose, without loss 
of generality, that the regular part of $\Sigma$ does not contain the 
origin; see Lemma~4.3.2 in \cite{Ru}. Let $\xi\neq0$ be any fixed point 
in the regular part of $\Sigma$. We may also suppose by simplicity that 
the first entry $\xi_1\neq0$, and so we define the following mapping 
$\eta:\cc^n\to\cc^n$ and variety $Y$,
\begin{equation}
\label{eqn1}\begin{array}{rcl}
\eta(y)&:=&(y_1/\xi_1)^\beta*(\xi_1,y_2,y_3,...,y_n),
\quad\hbox{for}\quad{y}\in\cc^n,\\
Y&:=&\{\widehat{y}\in\cc^{n-1}:Q_k(\xi_1,\widehat{y})=0,\,\forall\,k\}.
\end{array}\end{equation}

The action $s^{\beta}*z$ was given in (\ref{action}). We have that 
$\eta(\xi)=\xi$, and that the following identities hold for all $s\in\cc$ 
and $\widehat{y}\in\cc^{n-1}$ (recall equation~(\ref{homo}) and the fact 
that $\Sigma$ is the zero locus of the polynomials $\{Q_k\}$):
\begin{equation}
\label{eqn2}\begin{array}{rcl}
Q_k(\eta(s,\widehat{y}))&=&(s/\xi_1)^{d_k}\,
Q_k(\xi_1,\widehat{y}),\quad\hbox{and~so}\\
\eta(\cc^*\times{Y})&=&\{z\in\Sigma:z_1\neq0\}.
\end{array}\end{equation}

The symbol $\cc^*$ stands for $\cc\setminus\{0\}$. The mapping $\eta(y)$ 
is locally a biholomorphism whenever the first entry $y_1\neq0$. Whence, 
the point $\xi$ lies in the regular part of the variety $\cc\times{Y}$, 
because $\xi=\eta(\xi)$ also lies in the regular part of $\Sigma$ and 
$\xi_1\neq0$.  Thus, we can find a biholomorphism $\pi$ defined from an 
open domain $U$ in $\cc^m$ onto an open set in the regular part of $Y$, 
such that $\pi(\rho)$ is equal to $(\xi_2,...,\xi_n)$ for some 
$\rho\in{U}$. Consider the following holomorphic mapping defined 
for all points $s\in\cc$ and $x\in{U}$,
\begin{equation}\label{eqn3}
\Pi(s,x)\,:=\,s^{\beta}*(\xi_1,\pi(x))
\,=\,\eta(s\xi_1,\pi(x))\,\in\,\Sigma.
\end{equation}

The image $\Pi(\cc\times U)$ will be known as a {\bf generalized cone} 
from now on. Notice that $\Pi(\cc^*\times{U})$ lies in the regular part 
of $\Sigma$, for $\pi(U)$ is contained in the regular part of $Y$. The 
mapping $\Pi(s,x)$ is locally a biholomorphism whenever $s\neq0$, because 
$\eta$ is also a local biholomorphism for $y_1\neq0$. Finally, the image 
$\Pi(1,\rho)$ is equal to $\xi$.  Hence, recalling the differential form 
$\lambda$ and the function $g$ defined in~(\ref{fng-1}), we only need 
to prove that the pull-back $\Pi^*\lambda$ is equal to $\dbar{g}(\Pi)$ 
inside $\cc\times{U}$, in order to conclude that the $\dbar$-equation 
$\lambda=\dbar{g}$ holds in a neighborhood of $\xi$. Consider the 
following identity obtained by applying (\ref{action}) and (\ref{eqn3}) 
into (\ref{fng-1}), we define $\pi_1(x)\equiv\xi_1$,
\begin{equation}\label{eqn4}
g(\Pi(s,x))=\sum_{k=1}^n\frac{\beta_k}{2\pi{i}}\int_{\cc}
f_k(\Pi(ws,x))\frac{\overline{(ws)^{\beta_k}\pi_k(x)}\,
dw\wedge{d}\overline{w}}{\overline{w}\;(w-1)}.
\end{equation}
re
The given hypotheses on $\lambda$ yield that the pull-back $\Pi^*\lambda$ 
is $\dbar$-closed and bounded in $\cc^*\times{U}$, and so it is also 
$\dbar$-closed in $\cc\times{U}$; see Lemma~4.3.2 in \cite{Ru} or 
Lemma~(2.2) in \cite{SZ}. We can then use equations~(\ref{action}) 
and~(\ref{eqn3}) in order to calculate $\Pi^*\lambda$ when $\lambda$ 
is given by $\sum_kf_kd\overline{z_k}$,
\begin{eqnarray}
\nonumber\Pi^*\lambda&=&F_0(s,x)d\overline{s}\,+\,
\sum_{j\geq1}F_j(s,x)d\overline{x_j},\quad\hbox{with}\\
\label{eqn5}F_0(s,x)&=&\sum_{k=1}^nf_k(\Pi(s,x))\,
\beta_k\,\overline{s^{\beta_k-1}\pi_k(x)}.
\end{eqnarray}

Recall that $\pi_1(x)\equiv\xi_1$. Equation~(\ref{eqn3}) and the fact 
that $\lambda$ has compact support on $\Sigma$ also imply that the 
previous function $F_0(s,x)$ has compact support on every complex line 
$\cc\times\{x\}$, for all $x\in U$. Whence, the following Cauchy-Pompeiu 
integral satisfies the $\dbar$-equation $\Pi^*\lambda=\dbar{G}$ in the 
product $\cc\times U$, given $s\in\cc$ and $x\in U$,
\begin{equation}\label{eqn6}
G(s,x)\,:=\,\frac{1}{2\pi{i}}\int_{u\in\cc}
\frac{F_0(u,x)}{u-s}du\wedge{d}\overline{u}.
\end{equation}

Finally, equations~(\ref{eqn4}) and~(\ref{eqn6}) are identical, for we 
only need to apply the change of variables $u=sw$. Thus, the differential 
$\dbar{g}(\Pi)$ (resp.~$\dbar{g}$) is equal to the form $\Pi^*\lambda$ 
(resp.~$\lambda$) inside the space $\cc\times{U}$ (resp.~an open 
neighborhood of $\xi$); and so, the $\dbar$-equation $\lambda=\dbar{g}$ 
holds in the regular part of $\Sigma$, because $\xi\neq0$ was chosen 
in an arbitrary way in the regular part of $\Sigma$, and Lemma~4.3.2 
in \cite{Ru}.

\section{H\"older Estimates}\label{section:hoelder}

In this section, we will prove anisotropic H\"older estimates on the 
subvariety $\Sigma\subset\cc^n$ in the particular case when $\Sigma$ is 
homogeneous (a cone) and has got only one isolated singularity at the 
origin (Lemma~\ref{isotropic}). These estimates easily lead to optimal 
H\"older estimates on such varieties (Theorem~\ref{hoelder}). We will 
show later (in section \ref{final}) how we can use previous results in 
order to deduce H\"older estimates on weighted homogeneous varieties 
with an isolated singularity as well. The given hypotheses imply that 
$\Sigma\setminus\{0\}$ is a regular complex manifold in $\cc^n$. Consider 
the compact \textbf{link} $K$ obtained by intersecting $\Sigma$ with the 
unit sphere $bB$ of radius $\sqrt{n}$ and center at the origin in $\cc^n$. 
Notice that every point $\xi\in{K}$ has got at least one coordinate with 
absolute value $|\xi_k|\geq1$. We follow the proof of Theorem~\ref{main}.

Thus, given any point $\xi\in{K}$, we construct a generalized cone which 
contains it. For example, if the first entry $|\xi_1|\geq1$, we build the 
subvariety $Y_\xi$ as in~(\ref{eqn1}). Then, we consider a biholomorphism 
$\pi_\xi$ defined from an open set $U_\xi\subset\cc^m$ into a neighborhood 
of $(\xi_2,...,\xi_n)$ in $Y_\xi$, and the mapping $\Pi_\xi$ defined as 
in~(\ref{eqn3}) from $\cc\times{U_\xi}$ into $\Sigma$. We also restrict 
the domain of $\Pi_\xi$ to a smaller set $\cc\times{U''_\xi}$, where: 
$U''_\xi\Subset{U'_\xi}\Subset{U_\xi}$, the open set $U'_\xi$ is 
smoothly bounded, and $\pi_\xi(U''_\xi)$ is a convex open neighborhood 
of $(\xi_2,...,\xi_n)$ in $Y_\xi$. The generalized cone 
$\Pi_\xi(\cc{\times}U''_\xi)$ obviously contains to $\xi$, as we wanted. 
Recall that an open set $V$ in $Y_\xi$ is called convex whenever every 
pair of points in $V$ can be joined by a geodesic which is also contained 
in $V$. We proceed in a similar way when any other entry $|\xi_k|\geq1$.

Now then, since the link $K$ is compact, we may choose finitely 
many (let us say $N$) points $\xi^1,...,\xi^N$ in $K$ such that 
$K$ itself is covered by their associated generalized cones 
$C_j:=\Pi_{\xi^j}(\cc{\times}U''_{\xi^j})$. We assert that the 
analytic set $\Sigma$ is covered by the cones $C_j$. Let $z$ be any 
point in $\Sigma\setminus\{0\}$. It is easy to deduce the existence 
of $s\in\cc^*$ such that $s^\beta*z$ lies in $K$; and so there exists 
an index $1\leq{j}\leq{N}$ such that $s^\beta*z$ also lies in $C_j$. We 
may suppose that the first entry $|\xi^j_1|\geq1$, and that $\Pi_{\xi^j}$ 
is given as in~(\ref{eqn3}). Hence, there is a pair $(t,x)$ in the 
Cartesian product $\cc^*{\times}U''_{\xi^j}$ with
\begin{eqnarray*}
s^\beta*z&=&\Pi_{\xi^j}(t,x)\;=\;
t^\beta*(\xi^j_1,\pi_{\xi^j}(x));\quad\hbox{and~so}\\
z&=&(t/s)^\beta*(\xi^j_1,\pi_{\xi^j}(x))\;=\;\Pi_{\xi^j}(t/s,x).
\end{eqnarray*}

Previous identity shows that the whole analytic set $\Sigma$ is 
covered by the $N$ generalized cones $C_1, ..., C_N$. On the other hand, 
in order to prove the H\"older continuity of~(\ref{hoelder2}), we take 
a fixed parameter $0<\theta<1$ and a pair of points $z$ and $w$ in the 
intersection of $\Sigma$ with the open ball $B_R$ of radius $R>0$ and 
center at the origin in $\cc^n$. We want to show that there is a constant 
$C_\Sigma(\theta)>0$ which does not dependent on $z$ or $w$ such that:
\begin{equation}\label{hoelder5}
|g(z)-g(w)|\leq C_\Sigma(\theta)\cdot 
\dist_\Sigma(z,w)^\theta\cdot\|\lambda\|_\infty.
\end{equation}

One first step is to show that we only need to verify previous 
H\"older inequality when the points $z$ and $w$ are both contained in 
$B_R\cap{C_j}$, where $C_j$ is a unique generalized cone defined as in the 
paragraphs above. Let $\epsilon>0$ be a given parameter. The definition of 
$\dist_\Sigma(z,w)$ implies the existence of a piecewise smooth curve 
$\gamma_\epsilon:[0,1]{\rightarrow}\Sigma$ joining $z$ and $w$, i.e. 
$\gamma_\epsilon(0)=z$ and $\gamma_\epsilon(1)=w$, such that:
$$\hbox{length}(\gamma_\epsilon)=\int_0^1\|\gamma'(t)\|
\,dt \leq\dist_\Sigma(z,w)+\epsilon.$$

The image of $\gamma_\epsilon$ is completely contained in $B_R\cap\Sigma$ 
because $\Sigma$ is homogeneous (a cone). Now then, we are done if the 
points $z$ and $w$ are both contained in the same generalized cone $C_j$. 
Otherwise, we run over the curve $\gamma_\epsilon$ from $z$ to $w$, and 
pick up a finite set $\{z_k\}$ inside $\gamma_\epsilon\subset{B_R}$ such 
that: the initial point $z_0=z$, the final point $z_N=w$, two consecutive 
elements $z_j$ and $z_{j+1}$ lie in the same generalized cone, and three 
arbitrary elements of $\{z_k\}$ cannot lie in the same generalized cone. 
In particular, we may also suppose, without loss of generality, that: 
$z_0=z$ is in $C_1$, the final point $z_N=w$ is in $C_N$, and any other 
point $z_j$ is in the intersection $C_j\cap{C}_{j+1}$ for every index 
$1\leq{j}<N$. So that, two consecutive points $z_{j-1}$ and $z_j$ lie in 
the same generalized cone $C_j\cap{B_R}$ for each index $1\leq{j}\leq{N}$. 
Assume for the moment that there exist constants $C^j_\Sigma(\theta)>0$ 
such that
$$|g(z_{j-1}){-}g(z_j)|\leq{C^j_\Sigma}(\theta)\cdot 
\dist_\Sigma(z_{j-1},z_j)^\theta\cdot\|\lambda\|_\infty,$$
for all $1\leq{j}\leq{N}$. Then, it follows that
\begin{eqnarray*}
&&|g(z){-}g(w)|\leq\sum_{j=1}^N |g(z_{j-1}){-}g(z_j)|\leq\sum_{j=1}^N
C^j_\Sigma(\theta)\dist_\Sigma(z_{j-1},z_j)^\theta\|\lambda\|_\infty\\
&&\hspace{13ex}\leq{}C_\Sigma(\theta)\cdot
[\dist_\Sigma(z,w)+\epsilon]^\theta\cdot\|\lambda\|_\infty,
\end{eqnarray*}
where we have chosen $C_\Sigma(\theta)=\sum_jC^j_\Sigma(\theta)$. Since 
previous inequality holds for all $\epsilon>0$, it follows that we only 
need to prove that the H\"older estimates~(\ref{hoelder5}) holds under 
the assumption that $z$ and $w$ are both contained the intersection of 
a unique generalized cone $C_j$ with the open ball $B_R$ of radius $R>0$ 
and center at the origin in $\cc^n$. Moreover, we can suppose, without 
loss of generality, that $C_j$ is indeed the generalized cone given 
in~(\ref{eqn3}).

Recall the given hypotheses: The subvariety $\Sigma$ is homogeneous (a 
cone) and has got only one isolated singularity at the origin of $\cc^n$, 
so that each entry $\beta_k=1$ in definition~\ref{whs}. We fix a point 
$\xi$ in the link $K\subset\Sigma$, and assume that its first entry 
$|\xi_1|\geq1$. The subvariety $Y$ is then given in~(\ref{eqn1}), and 
the biholomorphism $\pi$ is defined from an open set $U\subset\cc^m$ 
into a neighborhood of $(\xi_2,...,\xi_n)$ in $Y$. Let $\lambda$ be a 
$(0{,}1)$-form as in the hypotheses of Theorem~\ref{main}. We may easily 
calculate the pull-back $\Pi^*\lambda$, with the mapping $\Pi$ given 
in~(\ref{eqn3}) for all $s\in\cc$ and $x\in{U}$,
\begin{equation}\label{hoelder6}
\Pi(s,x)=s^{(1,...,1)}*(\xi_1,\pi(x))=(s\xi_1,s\pi(x))\in\Sigma.
\end{equation}

The pull-back 
$\Pi^*\lambda=F_0(s,x)d\overline{s}+\sum_jF_jd\overline{x_j}$ satisfies:
\begin{eqnarray*}
F_0(s,x)&=&\sum_{k=1}^nf_k(\Pi(s,x))
\overline{\pi_k(x)},\quad\pi_1(x)\equiv\xi_1,\\
F_j(s,x)&=&\sum_{k=2}^nf_k(\Pi(s,x))\overline{
\Big[s\,\frac{\partial\pi_k}{\partial{x_j}}\Big]}.
\end{eqnarray*}

The hypotheses of Theorem~\ref{hoelder} yield that the support of 
every $f_k$ is contain in a ball of radius $R>0$ and center at the 
origin. Whence, equation~(\ref{hoelder6}) and the fact that $|\xi_1|\geq1$ 
automatically imply that each function $F_k(s,x)$ vanishes whenever 
$|s|>R$. Now then, we restrict the domain of $\Pi$ to a smaller set 
$\cc\times{U''}$, where: $U''\Subset{U'}\Subset{U}$, the open set 
$U'$ is smoothly bounded, and $\pi(U'')$ is a convex neighborhood 
of $(\xi_2,...,\xi_n)$ in $Y$.

Equation~(\ref{hoelder6}) yields that $\Pi_\xi(\{s\}{\times}U''_\xi)$ is 
convex for all $s\in\cc$ as well. The biholomorphism $\pi$ has also got a 
Jacobian (determinant) which is bounded from above and below (away from 
zero) in the compact closure $\overline{U'_\xi}$. Whence, there exists a 
constant $D_1>0$, such that the following identities hold for every point 
$(s,x)$ in $\cc\times\overline{U'}$ and each index $1\leq{j}\leq{m}$,
\begin{eqnarray}\label{eq:F0}
|F_0(s,x)|&\leq&D_1\cdot\|\lambda\|_\infty,\\
\label{eq:Fj}|F_j(s,x)|&\leq&D_1\cdot|s|\cdot\|\lambda\|_\infty.
\end{eqnarray}

We may show that the H\"older estimate~(\ref{hoelder5}) holds for 
all points $z$ and $w$ in the intersection of the generalized cone 
$\Pi(\cc{\times}U'')$ with the ball $B_R$, and so being able to conclude 
that the same estimate holds on $B_R\cap\Sigma$. Fix the parameter 
$0<\theta<1$. We are going to analyze two different cases. Firstly, we 
assume there exist a point $x\in{U''}$ and two complex numbers $s$ and 
$s'$, such that $z=\Pi(s,x)$ and $w=\Pi(s',x)$. We say, in this case, 
that $z$ and $w$ lie in the same complex line. Equation~(\ref{hoelder6}) 
and the fact that $|\xi_1|\geq1$ yields that $|s|$ is bounded:
\begin{equation}\label{hoelder7}
|s|\,\leq\,|s\xi_1|\,\leq\,\|z\|<R.
\end{equation}

The function $G=g(\Pi)$ defined 
in~(\ref{eqn4}) and~(\ref{eqn6}) satisfies:
\begin{eqnarray*}
&&|g(z)-g(w)|\;=\;|G(s,x)-G(s',x)|\\
&&\qquad=\;\frac{1}{2\pi}\bigg|\int_{|u|\leq{R}}F_0(u,x)\Big(
\frac{1}{u-s}-\frac{1}{u-s'}\Big)du\wedge{d}\overline{u}\bigg|.
\end{eqnarray*}

Recall that $F_0(u,x)$ vanishes whenever $|u|>R$. It is well known that 
there exists a constant $D_2(R,\theta)>0$, depending only on the radius 
$R>0$ and the parameter $\theta$, such that:
\begin{eqnarray}\label{eq:h1}
|g(z)-g(w)|\,\leq\,D_2(R,\theta)\,
|s-s'|^{\theta}\,D_1\,\|\lambda\|_\infty.
\end{eqnarray}

Notice that we have used~(\ref{eq:F0}), and consider chapter~6.1 of 
\cite{Ru} for a (more general) version of the inequality above. The 
analysis done in the previous paragraphs shows that~(\ref{hoelder5}) 
holds in the first case. Besides, since both $g(0)$ and $\Pi(0,x)$ 
vanishes, we also obtain the following useful estimate:
\begin{eqnarray}\label{eq:h2}
|G(s,x)|\,=\,|g(z)|\,\leq\,D_2(R,\theta)
\,D_1\,|s|^{\theta}\,\|\lambda\|_\infty.
\end{eqnarray}

We analyze now the symmetrical case. Let $z$ and $\widehat{w}$ be a 
pair of points in the intersection of $\Pi(\cc{\times}U'')$ with the 
ball $B_R$. Assume there exist a complex number $s\neq0$ and a pair of 
points $x$ and $x'$ in the open set $U''$ such that $z=\Pi(s,x)$ and 
$\widehat{w}=\Pi(s,x')$. We say, in this case, that $z$ and $\widehat{w}$ 
lie in the same slice. By a unitary change of coordinates which does not 
spoil the inequality~\eqref{eq:Fj}, we may assume that the entries of $x$ 
and $x'$ are all equal, with the possible exception of the first one. 
That is, we may assume that both $x$ and $x'$ lie in the complex line 
$L:=\cc\times\{(x_2,...,x_m)\}$. Recall that the differential $\dbar{G}$ 
is equal to $\Pi^*\lambda$ in the open set $\cc\times{U}$, according to 
equation~(\ref{eqn6}) and the statement just above it. Hence, we can 
evaluate $g(z)$ via the inhomogeneous Cauchy-Pompeiu formula on the 
line $L$,
\begin{eqnarray*}
g(z)=G(s,x)=\frac{1}{2\pi{i}}\int_{L\cap{U'}}
F_1(s,t,x_2,...,x_m)\frac{dt\wedge{d}\overline{t}}{t-x_1}&&\\
+\;\frac{1}{2\pi{i}}\int_{L\cap{bU'}}G(s,t,x_2, ..., x_m)
\frac{dt}{t-x_1},&&
\end{eqnarray*}
because $x$ is in $L\cap{U''}$ and $U''\Subset{U'}$. We introduce some 
notation in order to simplify the analysis. The symbols $I_1(s,x)$ and 
$I_2(s,x)$ stands for the above integrals on the set $L\cap{U'}$ and 
the boundary $L\cap{bU'}$, respectively. In particular, we have that:
$$g(\widehat{w})\,=\,G(s,x')\,=\,I_1(s,x')\,+\,I_2(s,x').$$

Recall that $x$ and $x'$ are both in $L\cap{U''}$, and that the 
difference $x-x'$ is equal to the vector $(x_1{-}x_1',0,...,0)$. 
Inequality \eqref{eq:Fj} implies the existence of a constant 
$D_3(\theta)>0$, depending only on the diameter of $U'$ and the 
parameter $\theta$, such that:
\begin{equation}\label{eq:h3}
|I_1(s,x)-I_1(s,x')|\,\leq\,D_3(\theta)\,
|x_1-x_1'|^\theta\,D_1\,|s|\,\|\lambda\|_\infty.
\end{equation}

We can calculate similar estimates for $I_2$. Let $\delta>0$ be the 
distance between the compact sets $\overline{U''}$ and $bU'$ in $\cc^m$. 
We obviously have that $\delta>0$ because $U''\Subset{U'}$. The following 
estimates are deduced from~(\ref{eq:h2}) and the mean value Theorem, the 
maximum is calculated over all $u$ in $L\cap{U'}$,
\begin{eqnarray*}
|I_2(s,x){-}I_2(s,x')|\leq\frac{|x_1{-}x_1'|}{2\pi}
\max_{v}\bigg|\int_{L\cap{bU'}}\hspace{-1ex}
\frac{G(s,t,x_2,...)dt}{(t-v)^2}\bigg|&&\\
\leq\frac{|x_1{-}x_1'|}{2\pi}\cdot\frac{\hbox{length}(L\cap{bU'})}
{\delta^2}\,D_2(R,\theta)\,D_1\,|s|^\theta\,\|\lambda\|_\infty.&&
\end{eqnarray*}

Previous estimates and the inequalities (\ref{eq:h1}) and (\ref{eq:h3}) 
can be summarized in the following lemma. It is convenient to recall 
that the points $x$ and $x'$ are both contained in the bounded set 
$U''\Subset\cc^m$. Moreover, we also have that $|s|<R$ and $|s'|<R$, 
because $z$, $w$ and $\widehat{w}$ are all contained in the ball $B_R$; 
recall the proof of~(\ref{hoelder7}).

\begin{lemma}[\textbf{Isotropic Estimates}]\label{isotropic} 
In the situation of Theorems~\ref{main} and~\ref{hoelder}, consider the 
functions $g$ and $\Pi$ given in~(\ref{fng-1}) and~(\ref{hoelder6}), 
respectively, and the bounded open set $U''\Subset\cc^m$ defined in 
the paragraphs above. Then, for every parameter $0<\theta<1$, there is 
a constant $D_4(R,\theta)>0$ which does not depend on $\lambda$ such 
that the following statements hold for all the points $z=\Pi(s,x)$ and 
$w=\Pi(s',x')$ in the intersection of $\Pi(\cc{\times}U'')$ with the 
ball $B_R$:
$$|g(z)-g(w)|\,\leq\,D_4(R,\theta)\,|s-s'|^\theta\,\|\lambda\|_\infty,$$
whenever $x=x'$, i.e. $z$ and $w$ are in the same line; and:
$$|g(z)-g(w)|\,\leq\,D_4(R,\theta)\,\|x-x'\|^\theta
\,|s|^\theta\,\|\lambda\|_\infty,$$
whenever $s=s'$, i.e. $z$ and $w$ are in the same slice.
\end{lemma}

It is now easy to prove that the H\"older estimates given 
in~(\ref{hoelder2}) and~(\ref{hoelder5}) hold for all of points $z$ and 
$w$ which fulfill the assumptions of Lemma~\ref{isotropic}; so that they 
lie in the intersection of the generalized cone $\Pi(\cc{\times}U'')$ 
with the ball $B_R$. The definition of $\Pi$, given in~(\ref{hoelder6}), 
allows us to write down the identities:
\begin{equation}\label{eq:h5}
z=\Pi(s,x)=s\,(\xi_1,\pi(x))\quad\hbox{and}
\quad{}w=\Pi(s',x')=s'\,(\xi_1,\pi(x')). 
\end{equation}

Fix the point $z':=\Pi(s,x')=s\,(\xi_1,\pi(x'))$ such that it is in 
the same line than $w$ and in same the slice than $z$. We can suppose, 
without loss of generality, that $z'\in{B_R}$ because $z$ and $w$ also 
lie in $B_R$. Otherwise, if the norm $\|z'\|\geq{}R$, we only need to 
use $\Pi(s',x)$ instead. We can easily deduce the following estimate 
from~(\ref{eq:h5}) and the fact that $|\xi_1|\geq1$,
$$|s-s'|\leq|s\xi_1-s'\xi_1|\leq\|z-w\|\leq\dist_\Sigma(z,w).$$

Recall that $\pi$ is a biholomorphism whose Jacobian (determinant) 
is boun\-ded from above and below (away from zero) in the compact set 
$\overline{U''}$. Moreover, the image $\pi(U'')$ is also a convex set 
in $Y$. Whence, recalling~(\ref{eq:h5}), we can deduce the existence 
of a constant $D_5>0$, depending only on $\pi$ and $U''$, such that:
\begin{eqnarray*}
&&\frac{|s|\cdot\|x{-}x'\|}{D_5}\leq|s|
\cdot\|\pi(x){-}\pi(x')\|\leq\|z{-}z'\|\\
&&\leq\|z{-}w\|+\|w{-}z'\|\leq\|z{-}w\|
+|s{-}s'|\cdot\|(\xi_1,\pi(x'))\|\\
&&\leq\dist_\Sigma(z,w)\cdot\big[2+\|\pi(x')\|\big].
\end{eqnarray*}

Thus, there exists a constant $D_6>0$, depending only on $\pi$ and $U''$, 
such that the following identities hold for all the points $z=\Pi(s,x)$ 
and $w=\Pi(s',x')$ in the intersection of $\Pi(\cc{\times}U'')$ with 
the ball $B_R$,
$$|s{-}s'|\leq{D_6}\cdot\dist_\Sigma(z,w)\quad\mbox{and}
\quad|s|\cdot\|x{-}x'\|\leq{D_6}\cdot\dist_\Sigma(z,w).$$

Recall that $z'$ is in the same line than $w$ and in same the slice 
than $z$; so that Lemma \ref{isotropic} automatically yields that:
\begin{eqnarray*}
&&|g(z)-g(w)|\,\leq\,|g(z)-g(z')|+|g(z')-g(w)|\\
&&\leq\,D_4(R,\theta)\Big[|s|^\theta \cdot\|x-x'\|^\theta
+|s-s'|^\theta\Big]\|\lambda\|_\infty\\
&&\leq\,D_4(R,\theta)\,D_6^\theta\,
\dist_\Sigma(z,w)^\theta\,\|\lambda\|_\infty.
\end{eqnarray*}

This completes the proof that the H\"older estimates given 
in~(\ref{hoelder2}) and~(\ref{hoelder5}) hold for all of points $z$ 
and $w$ in the intersection of the ball $B_R$ with the generalized cone 
$\Pi(\cc{\times}U'')$; and so we can conclude that the same H\"older 
estimates hold for all point $z$ and $w$ in $B_R\cap\Sigma$, we just 
need to recall the analysis done in the paragraphs located between 
equations~(\ref{hoelder5}) and~(\ref{hoelder6}).

\section{$L^2$-Estimates}\label{section:L2}

We prove Theorem~\ref{thm:l2} in this section, so we begin by showing 
that the function $g$ given in (\ref{eq:l2}) is indeed well defined. 
Recall that $\Sigma$ is a pure $d$-dimensional homogeneous (cone) 
subvariety of $\cc^n$, so that $n\geq2$ and each entry $\beta_k=1$ in 
Definition~\ref{whs}.  Besides, consider the differential form $\lambda$ 
given by $\sum_kf_kd\overline{z_k}$, where the coefficients $f_k$ are 
all square-integrable functions in $\Sigma$. Assume that the support of 
$\lambda$ is contained in the open ball $B_R$ of radius $R>0$ and center 
at the origin. We only need to show that the following integrals exist,
\begin{equation}\label{eqn10}
\int_{w\in\cc}\int_{z\in\Sigma\cap{B_R}}\Big|\frac{f_k
(wz)w^dz_k}{w(w-1)}\Big|\,dV_\Sigma\,dV_\cc\,<\,\infty.
\end{equation}

A direct application of Fubini's theorem will yield that the 
integrals in (\ref{eq:l2}) are all well defined for almost all $z$ in 
$\Sigma\cap{B_R}$, and so they are also well defined for almost all 
$z\in\Sigma$ because the radius $R$ can be as large as we want. The fact 
that $\Sigma$ is a pure $2d$-real dimensional and homogeneous (cone) 
subvariety automatically implies the existence of a constant $C_0>0$ such 
that the following equations hold for all $w\in\cc$ and real $\rho>0$:
\begin{equation}\label{eqn11}
\int_{z\in\Sigma\cap{B_\rho}}\!\|z\|^2=C_0\rho^{2d+2},\quad
\int_{z\in\Sigma}|f_k(wz)|^2=\frac{\|f_k\|^2_{L^2(\Sigma)}}{|w|^{2d}}.
\end{equation}

Notice that we need not calculate the integral (\ref{eqn10}) in the 
Cartesian product of $\cc$ times $\Sigma\cap{B_R}$. We can simplify the 
calculations by integrating over the set $\Xi$ defined below, because 
$f_k(wz)=0$ whenever $\|wz\|\geq{R}$,
\begin{equation}\label{eqn12}
\Xi\,:=\,\{(w,z)\in\cc\times\Sigma\,:\,\|z\|<R,\,\|wz\|<R\}.
\end{equation}

We easily have that:
\begin{eqnarray}\label{eqn13}
\Big\|\frac{f_k(wz)w^d}{|w^2-w|^{2/3}}\Big\|_{L^2(\Xi)}^2\;\leq\;
\int_{w\in\cc}\int_{z\in\Sigma}\frac{|f_k(wz)w^d|^2}{|w^2-w|^{4/3}}&&\\
\nonumber\leq\;\|f_k\|^2_{L^2(\Sigma)}\int_{w\in\cc}\frac{1}
{|w^2-w|^{4/3}}\;<\;\infty;&&
\end{eqnarray}
and that:
\begin{eqnarray}\label{eqn14}
&&\Big\|\frac{z_k}{|w^2-w|^{1/3}}\Big\|^2_{L^2(\Xi)}\;\leq\;
\int_{w\in\cc}\int_{\genfrac..{0pt}1{z\in\Sigma\cap{B_R},}
{\|z\|<R/|w|}}\frac{\|z\|^2}{|w^2-w|^{2/3}}\\
\nonumber&&\leq\;\int_{|w|\leq1}\frac{C_0R^{2d+2}}{|w^2-w|^{2/3}}
\;+\;\int_{|w|>1}\frac{C_0(R/|w|)^{2d+2}}{|w^2-w|^{2/3}}\;<\;\infty.
\end{eqnarray}

The last integral in the first line of (\ref{eqn14}) must be separated 
into two parts according to the fact that $|w|$ is either less or greater 
than one, and then, one must apply (\ref{eqn11}) with $\rho$ respectively 
equal to $R$ or $R/|w|$. The Cauchy-Schwartz inequality 
$\|ab\|_{L^1}\leq\|a\|_{L^2}\|b\|_{L^2}$ allows us to deduce equation 
(\ref{eqn10}) from the inequalities (\ref{eqn13})--(\ref{eqn14}), for we
only need to integrate on the set $\Xi$ given in~(\ref{eqn12}). Thus, a 
direct application of Fubini's theorem implies that the following function 
given in~(\ref{eq:l2}) is well defined for almost all $z\in\Sigma$,
\begin{equation}\label{eqn15}
g(z)=\sum_{k=1}^n\frac{H_k(z)}{\pi},\quad
H_k(z)=\int_{|w|\leq\frac{R}{\|z\|}}\!f_k(wz)
\frac{w^d\,\overline{z_k}\,dw\wedge{d}\overline{w}}{w(w-1)\,2i}.
\end{equation}

Notice that $f_k(wz)=0$ whenever $\|wz\|\geq{R}$, because of the given 
hypotheses. Moreover, the $L^2$-estimate~(\ref{eq:l3}) easily follows 
from the following inequalities (Cauchy-Schwartz) and equations 
(\ref{eqn13})--(\ref{eqn14}):
\begin{eqnarray*}
&&|H_k(z)|^2\leq\int_{|w|\leq\frac{R}{\|z\|}}\frac{|f_k(wz)w^d|^2}
{|w^2-w|^{4/3}}\int_{|\omega|\leq\frac{R}{\|z\|}}\frac{|z_k|^2}
{|\omega^2-\omega|^{2/3}}\quad\hbox{and}\\
&&\int_{z\in\Sigma\cap{B_R}}\!|H_k(z)|^2\leq
\Big\|\frac{f_k(wz)w^d}{|w^2-w|^{2/3}}\Big\|_{L^2(\Xi)}^2
\Big\|\frac{z_k}{|w^2-w|^{1/3}}\Big\|^2_{L^2(\Xi)}.
\end{eqnarray*}

Recall that $\big|\frac{dw{\wedge}d\overline{w}}{2i}\big|$ is the volume 
differential in $\cc$, that $\|f_k\|_{L^2(\Sigma)}$ is less than equal 
to $\|\lambda\|_{L^2_{0,1}(\Sigma)}$, and that we only need to integrate 
on the set $\Xi$ given in~(\ref{eqn12}). Finally, we prove that $g$ in 
(\ref{eqn15}) satisfies the differential equation $\dbar{g}=\lambda$. 
We only need to follow steep by the steep the proof presented in 
Section~\ref{section:main}. The only difference is that we must use a 
weighted Cauchy-Pompeiu integral in~(\ref{eqn6}), with $m=d{-}1$ integer:
\begin{eqnarray}\label{eqn16}
\mathcal{G}(s,x)&:=&\frac{1}{2\pi{i}}\cdot\frac{1}{s^m}
\int_{u\in\cc}\frac{u^mF_0(u,x)}{u-s}du\wedge{d}\overline{u},\\
\nonumber\hbox{where}&&F_0(u,x)\,=\,
\sum_{k=1}^nf_k(\Pi(u,x))\,\overline{\pi_k(x)}.
\end{eqnarray}

Notice that that $\Pi(u,x)=u(\xi_1,\pi(x))$ because each entry 
$\beta_k=1$ in~(\ref{action}) and~(\ref{eqn3}). We obviously have that 
$\dbar\mathcal{G}=[s^m\Pi^*\lambda]/s^m$. Hence, the function $g$ given in 
(\ref{eqn15}) is a solution to $\dbar{g}=\lambda$, because $g(\Pi(s,x))$ 
is identically equal to (\ref{eqn16}) after setting $u=sw$ and 
$\pi_1(x)\equiv\xi_1$. This concludes the proof of Theorem \ref{thm:l2}.

\section{Weighted Homogeneous Estimates}\label{final}

We want to close this paper presenting a useful technique for generalizing 
the estimates given in Theorems~\ref{hoelder} and~\ref{thm:l2}, so as to 
consider weighted homogeneous subvarieties instead of cones. Let 
$X\subset\cc^n$ be a weighted homogeneous subvariety with only one 
singularity at the origin and defined as the zero locus of a finite set 
of polynomials $\{Q_k\}$. Thus, the polynomials $Q_k(x)$ are all weighted 
homogeneous with respect to the same vector $\beta\in\zz^n$, and each 
entry $\beta_k\geq1$. Define the following holomorphic mapping:
\begin{equation}\label{eqn20}
\Theta:\cc^n\to\cc^n,\quad\hbox{with}\quad
\Theta(z)=(z_1^{\beta_1},z_2^{\beta_2},...,z_n^{\beta_n}).
\end{equation}

It is easy to see that each polynomial $Q_k(\Theta)$ is homogeneous, 
and so the subvariety $\Sigma\subset\cc^n$ defined as the zero locus 
of $\{Q_k(\Theta)\}$ is a cone. Moreover, since $\Theta$ is locally a 
biholomorphism in $\cc^n\setminus\{0\}$, we have that $\Sigma$ has got 
only one singularity at the origin as well. Consider a $(0{,}1)$-form 
$\aleph$ given by the sum $\sum_kf_kd\overline{x_k}$, where the 
coefficients $f_k$ are all Borel-measu\-rable functions in $X$, 
and $x_1,..., x_n$ are the Cartesian coordinates of $\cc^n$. We 
may follow two different paths in order to solve the equation 
$\overline{\partial}h=\aleph$. We may apply the main 
Theorem~\ref{main}, whenever $\aleph$ is bounded and 
has compact support on $X$, so as to get the solution:
$$h(x)=\sum_{k=1}^n\frac{\beta_k}{2\pi{i}}\int_{w\in\cc}
f_k(w^{\beta}*x)\frac{(\overline{w^{\beta_k}x_k})\,dw
\wedge{}d\overline{w}}{\overline{w}\;(w-1)}.$$

Otherwise, we may consider the pull-back $\Theta^*\aleph$, 
and apply Theorem~\ref{hoelder}, in order to solve the equation 
$\overline{\partial}g=\Theta^*\aleph$ on $\Sigma$. We easily have that:
\begin{eqnarray*}
\Theta^*\aleph&=&\sum_{k=1}^nf_k(\Theta(z))\,\beta_k
\overline{z_k}^{\beta_k-1}\,d\overline{z_k}\qquad\hbox{and}\\ 
g(z)&=&\sum_{k=1}^n\frac{\beta_k}{2\pi{i}}\int_{w\in\cc}
f_k(\Theta(wz))\frac{(\overline{wz_k})^{\beta_k}\,dw
\wedge{}d\overline{w}}{\overline{w}\;(w-1)}.
\end{eqnarray*}

Both paths yield exactly the same solution because $g(z)$ is identically 
equal to $h(\Theta(z))$. Recall that $w^{\beta}*\Theta(z)$ is equal to 
$\Theta(wz)$ for all $w\in\cc$ and $z\in\cc^n$. Hence, we may calculate 
the solution $g(z)$ above, and use the H\"older estimates given in 
equation~(\ref{hoelder2}),
$$|g(z)-g(w)|\leq{C}_\Sigma(R,\theta)\cdot 
\dist_\Sigma(z,w)^\theta\cdot\|\Theta^*\aleph\|_\infty.$$

A final steep is to \textit{push forward} these estimates, in order to 
deduce similar H\"older estimates for the solution $h(x)$ on $X$. A 
detailed analysis on the procedure for \textit{pushing forward} the H\"older 
estimates can be found in \cite{SZ}. On the other hand, we may use a similar 
procedure for $L^2$-estimates. In that case the subvariety $X$ can have arbitrary singularities.

\bibliographystyle{plain}

\end{document}